\documentclass{amsart}
\usepackage{url}

\newcommand{\<}{\langle}
\renewcommand{\>}{\rangle}
\newcommand{\ev}{\text{ev}}
\newcommand{\ad}{\text{ad}}

\renewcommand{\a}{\alpha}
\renewcommand{\b}{\beta}
\renewcommand{\c}{\gamma}
\renewcommand{\d}{\delta}

\newcommand{\Reals}{\mathbb R}
\renewcommand{\P}{{\mathcal P}}

\newcommand{\Z}{\mathbb Z}
\newcommand{\R}{\mathbb R}
\newcommand{\C}{\mathbb C}

\newcommand{\G}{\mathcal G}

\newcommand{\LK}{\mathfrak k}
\newcommand{\cstar}{\C^\times}
\newcommand{\Cstar}{\C^\times}

\theoremstyle{plain}
\newtheorem{theorem}{Theorem}[section]

\theoremstyle{definition}
\newtheorem{definition}[theorem]{Definition}

\theoremstyle{remark}

\begin{document}

\title[Bundle gerbes]{Higgs fields, bundle gerbes and string
structures}
\author{Michael K. Murray}
\address[Michael K. Murray]
{Department of Pure Mathematics\\
University of Adelaide\\
Adelaide, SA 5005 \\
Australia}
\email[Michael K. Murray]{mmurray@maths.adelaide.edu.au}
\author{Daniel Stevenson}
\email[Daniel Stevenson]{dstevens@maths.adelaide.edu.au}
\thanks{Both authors acknowledge the support of the Australian
Research Council.}

\subjclass{53C80, 81T50,
83E30}

\begin{abstract}
We use bundle gerbes and their connections and curvings to obtain
an explicit formula for a de Rham representative of the string
class of a loop group bundle. This is related to earlier work
on calorons.
\end{abstract}
\maketitle

\section{Introduction}
In this paper we bring together calorons (monopoles for the loop
group), bundle gerbes and string structures to produce a formula for
the string class of a principal bundle with structure group
the loop group.   If $K$ is a compact Lie group and $L (K)$ is
the group of smooth maps $\gamma$ from $[0, 2\pi]$ to $K$ such that
$\gamma(0) = \gamma(2\pi)$ it is well known \cite{PreSeg} that there 
is a central
extension
$$
0 \to U(1) \to \widehat{L(K)} \to L(K) \to 0
$$
where $\widehat{L(K)}$ is the Kac-Moody group.
If $P \to M$ is a principal bundle with structure group $L (K)$
it has a characteristic class, the {\em string class}, in $H^3(M,
\Z)$ which is the  obstruction to lifting $P$ to a $\widehat{L(K)}$
bundle.

In \cite{Mur1} Murray introduced the lifting bundle gerbe.
This is a bundle gerbe whose Dixmier-Douady class is the obstruction
to a principal $G$ bundle lifting to a principal $\hat G$ bundle
when
$$
0 \to \cstar \to \hat G \to G \to 0
$$
is a central extension.  In the case that $G = L (K)$ the
Dixmier-Douady class of the lifting bundle gerbe is the string class.
We use the methods of \cite{Mur1} to calculate
an explicit formula for the Dixmier-Douady class of the lifting
bundle gerbe when $G$ is the loop group and hence provide an
explicit differential three form
representative for the de Rham image of the string class in real cohomology.
This three form is defined in terms of a connection and
Higgs field for the loop group bundle.

In \cite{GarMur} it was shown that there is a
  bijective correspondence between  principal bundles over a
manifold $M$ with structure group  $L (K)$ and $K$ bundles over
$S^1 \times M$.
This was used to set up a correspondence between periodic instantons,
or  calorons, and loop group valued monopoles. In particular a
connection  for the $K$ bundle corresponded to a connection and Higgs
field for the $L(K)$ bundle.    We apply this correspondence
to show that the string class of an $L(K)$ bundle on $M$ is the
integral over the circle of
the Pontrjagin class of the
corresponding $K$ bundle over $S^1 \times M$.

Finally we relate these results
to earlier work  \cite{Kil, CarMur} on string structures.  Recall
that if
$Q$ is a  $K$ bundle  over a manifold $X$ we can take loops
everywhere and form a loop
group bundle $P= L (Q)$ over $M = L (X)$.  In that case it was known
from work of Killingback \cite{Kil}
that the Pontrjagin class of the bundle $Q$ in $H^4(X, \Z)$
transgressed to define the string class in $H^3(L(X), \Z)$.
The transgression consists of pulling back by the evaluation map
$$
\ev \colon S^1 \times L(X) \to X
$$
and pushing down to $L (X)$ by integrating over the circle.
If we apply the correspondence of \cite{GarMur} to the
principal $L (K)$ bundle $L (Q)(L (X), L (K))$ it produces a $K$
bundle over
$S^1 \times L (X)$ and
we show that this   is just the pull-back of $Q$ by the evaluation
map.
This recovers the result of Killingback \cite{Kil}.

As the present work was being completed a preprint was received from
Kiyonori Gomi \cite{Gom} which also defines connections
and curvings on lifting bundle gerbe  using the notion of {\em 
reduced splittings}
and building on results in \cite{Bry}.
We discuss the relationship between reduced splittings and Higgs fields in
\ref{sec:gomi}.

\section{Some preliminaries}

\subsection{$\cstar$ bundles}

Let $P \to X$ be a $\cstar$ bundle over a manifold $X$. We shall
denote
the fibre of $P$ over $x \in X$ by $P_x$.
Recall \cite{Bry} that if $P$ is a
  $\cstar$ bundle over a manifold $X$ we can define the dual
bundle $P^*$ as the same space $P$
but with the action $p^* g = (pg^{-1})^*$ and, that if $Q$ is another
such bundle, we can define the product bundle $P\otimes Q$ by
$(P\otimes Q)_x = (P_x \times Q_x)/\cstar $ where $\cstar$ acts
by $(p,q)w = (pw, qw^{-1})$. We denote an element of $P\otimes Q$ by
$p \otimes q$ with the understanding that $(pw) \otimes q = p \otimes
(qw) = (p\otimes q)w$ for $w \in \cstar$. It is straightforward to
check that
$P\otimes P^*$ is canonically trivialised by the section  $x \mapsto
p
\otimes p^*$ where $p$ is any point in $P_x.$

If $P$ and $Q$ are $\cstar$ bundles on $X$ with connections
$\mu_P$ and $\mu_Q$ then $P \otimes Q$ has an induced
connection we denote by $\mu_P \otimes \mu_Q$. The curvature
of this connection is $R_P + R_Q$ where $R_P$ and $R_Q$ are the
curvatures of $\mu_P$ and $\mu_Q$ respectively.
The bundle $P^*$ has an induced connection whose curvature
  is $-R_P$.

  \subsection{Simplicial spaces}
Recall \cite{Dup} that a  simplicial manifold $X$ is a collection of
spaces
$X_0, X_1, X_2, X_3, \dots$
  with maps $d_i \colon X_{p} \to X_{p-1}$ for $i = 1, \dots, p$, and
  $s_j \colon X_p \to X_{p+1}$ for $j=1, \dots, p$, satisfying
the, so-called, {\em simplicial identities}:
  \begin{align}
\label{eq:identities}
      d_i d_j &= d_{j-1}d_i ,  &i < j \\
      s_i s_j &= s_{j+1} s_i,  &i \leq j \\
     d_i s_j &= \begin{cases}
     s_{j-1}d_i,   &i < j \\
     \text{id}, &i = j, i=j+1\\
     s_j d_{i-1}, &i > j+1 .
     \end{cases}
     \end{align}

Let $\Omega^p(M)$ denote the space of all differentiable
$p$ forms on a manifold $M$. Define a homomorphism $\d\colon
\Omega^{n}(
X_p)\to \Omega^{n}(X_{p+1})$ by
$$
\d = \sum^{p} _{i=1}(-1)^{i-1}d_{i}^{*}.
$$
It is straightforward  to check that $\d^{2} = 0$ and it clearly
commutes with exterior derivative $d$. Hence we have a  complex
\begin{equation}
\label{eq:shortcomplex}
\Omega^{n}(X_0)\stackrel{\d}{\to} \Omega^{n}(X_1)
\stackrel{\d}{\to} \Omega^{n}(X_2)\stackrel{\d}
{\to}\cdots \stackrel{\d}{\to} \Omega^{n}(X_p)
\stackrel{\d}{\to} \cdots
\end{equation}

We remark that from a simplicial space we can
define a topological space called its {\it realisation}\footnote{It makes
no difference for the discussion in this paper if
it is the fat or geometric realisation.}. The double
complex $\Omega^p(X^q)$ with the differentials $d$ and $\delta$
has total cohomology the real cohomology of the realisation.

If $P \to X_p$ is a $\cstar$ bundle then we can define a $\cstar$
bundle over $X_{p+1}$ denoted  $\delta(P)$ by
$$
\delta(P) = d_1^{-1}(P) \otimes d_2^{-1}(P)^* \otimes
d_3^{-1}(P)
\otimes \dots.
$$
If $s$ is a section of $P$ then it defines $\delta(s)$ a section of
$\delta(P)$ and if $\mu$ is a connection on $P$ with curvature
$R$ it defines  a connection $\delta(\mu)$ on $\delta(P)$ with
curvature $\delta(R)$.
If we consider $\delta(\delta(P))$ it is a product of factors
and because of the simplicial identities
\eqref{eq:identities}  every factor occurs with its dual so 
$\delta(\delta(P))$ is
canonically
trivial.  If $s $ is a section of $P$ then under this identification
$\delta\delta(s) = 1$ and moreover if $\mu$ is a connection on $P$
then $\delta\delta(\mu) $ is the flat connection on
$\delta\delta(P)$ with respect to $\delta(\delta(s))$.

If $X$ is a simplicial space then  a {\em simplicial
line bundle} \cite{BryMcl} is a $\cstar$ bundle $P$ over $X_1$ with a section
$s$ of $\delta(P)$ over $X_2$ with the property that $\delta(s)$
is the canonical section of $\delta^2(P)$.

\subsection{Locally split maps}
We will be interested in maps  $\pi \colon Y \to M$ which admit local
sections.  That is, for every $x \in M$ there is an open set $U  $
containing $x$ and  a local section $s \colon U \to Y$.
  For want of a better term we will call maps like this
locally split.   Note that a locally split map is necessarily
surjective and that if we are dealing with the smooth
category a locally split map is just a submersion.
Locally trivial fibrations are, of course, locally split
but the converse is not true.

Let $Y \to M$ be locally split. Then we denote by
  $Y^{[2]} = Y\times_\pi Y$ the fibre product
of $Y$ with itself over $\pi$,
  that is the subset of pairs $(y, y')$ in
$Y \times Y$ such that $\pi(y) = \pi(y')$.  More generally we denote
the $p$th fold fibre product by $Y^{[p]}$.

For  $p = 1,2,\ldots$ we have $p$ projection maps
$\pi_{i}\colon Y^{[p+1]}\to Y^{[p]}$ for
$i=1,2, \ldots, p+1$ given by omitting the $i$-th
factor, so
$$
\pi_{i}(y_{1},y_{2},\ldots,y_{p+1})
= (y_{1},\ldots,y_{i-1},y_{i+1},\ldots,y_{p+1}).
$$
The spaces $ X_0 = Y, X_1 = Y^{[2]}, \dots$ define a simplicial
manifold.
For this simplicial manifold we can augment  the complex
\eqref{eq:shortcomplex} by adding at the beginning
the space $\Omega^n(M)$ and the map which
pulls back forms from $M$  to $Y$ to obtain a complex \begin{equation}
\label{eq:complex}
\Omega^{n}(M)\stackrel{\pi^{*}}{\to} \Omega^{n}(Y)
\stackrel{\d}{\to} \Omega^{n}(Y^{[2]})\stackrel{\d}
{\to}\cdots \stackrel{\d}{\to} \Omega^{n}(Y^{[p]})
\stackrel{\d}{\to} \cdots
\end{equation}

  It is a fundamental result of \cite{Mur1} that the  complex
\eqref{eq:complex},
  has no cohomology.
So if $\eta \in \Omega^{n}
(Y^{[p]})$ satisfies $\d(\eta) = 0$ then we can
solve the equation $\eta = \d(\rho)$ for some
$\rho \in \Omega^{n}(Y^{[p-1]})$ (we define $Y^{[0]}
= M$).  Note that a solution $\rho$ is not unique,
any two solutions $\rho$ and $\rho^{'}$ will differ
by $\d(\zeta)$ for some $\zeta\in \Omega^{n}(Y^{[p-1]})$.

\section{Central extensions}
We give here details of a method of constructing
central extensions presented in \cite{MurSte}.
Let $\G$ be a Lie group.
Recall that from $\G$  we can construct a
simplicial manifold $N\G = \{N\G_{p}\}$
with $N\G_{p} = \G^{p}$ with face operators $d_i \colon \G^{p+1}
\to \G^{p}$ defined by
$$
d_{i}(g_{1},\ldots,g_{p+1}) = \begin{cases}
                               (g_{2},\ldots,g_{p+1}), & i = 0, \\
                               (g_{1},\ldots,g_{i-1}g_{i},g_{i+1},
                                 \ldots,g_{p+1}), & 1\leq i\leq p-1,
\\
                                (g_{1},\ldots,g_{p}), & i = p.
                               \end{cases}
$$
Consider a central extension
$$
\cstar\to \hat{\G}\stackrel{\pi}{\to} \G.
$$
Following Brylinski and McLaughlin \cite{BryMcl}
we think of this as a $\cstar$ bundle $\hat\G \to \G$
with a product $M \colon \hat \G \times \hat  \G
\to \hat\G$ covering the product $m = d_1 \colon G \times G \to G$.

Because this is a central extension we must have that
$M(pz, qw) = M(p,q)zw$ for any $p, q \in P$ and $z, w \in
\cstar$.  This means we have a  section $s$
of $\delta(P)$ given by
$$
s(g, h)  = p \otimes M(p, q) \otimes q
$$
for
any $p \in P_g$ and $q \in P_h$.  This is well-defined as $pw \otimes
M(pw, qz) \otimes qz = pw \otimes
M(p, q)(wz)^{-1} \otimes qz = p \otimes  M(p, q) \otimes q$.
Conversely any such section gives rise to an $M$.

Of course we need an associative product and it can be
shown that $M$ being associative is equivalent to $\delta(s) = 1$.
To actually make $\hat \G$ into a group we need more
than multiplication we need an identity $\hat e \in \hat \G$ and
an inverse map.  It is straightforward to check that if $e \in \G$
is the identity then, because $M \colon \hat \G_e \times \hat \G_e
\to \hat\G_e$,  there is a unique $\hat e \in \hat \G_e$ such that
$M(\hat e, \hat e) = \hat e$.  It is also straightforward to
deduce the existence of a unique inverse.

Hence we have the result from \cite{BryMcl} that
a central extension of $\G$ is a $\cstar$ bundle $P \to G$
together with a section $s $ of $\delta(P) \to G\times G$
such that $\delta(s) = 1$.  In \cite{BryMcl} this is phrased
in terms of simplicial line bundles.

For our purposes we need to phrase this result in terms of
differential forms.
We call a connection for $\hat \G \to \G$, thought of as a
$\cstar$ bundle, a connection for the central extension.  An
isomorphism of central extensions with connection is an isomorphism
of bundles with connection which is a group isomorphism on the total
space $\hat\G$.  Denote by
$C(\G)$ the set of all isomorphism classes of
central extensions of $\G$ with connection.

Let  $\mu \in \Omega^1(\hat \G)$  be a connection on the
bundle  $\hat{\G} \to \G$ and consider the form
$$
\hat\delta(\mu) =  d_0^*(\mu) - d_1^*(\mu) + d_2^*(\mu)
$$
on $\hat\G \times \hat\G$.  There is a projection map
$$
\hat\G \times \hat\G \to \hat\G\otimes\hat\G
$$
whose fibre at $p \otimes q$ is all $(pz, qz^{-1})$ for
$z\in \Cstar$. If $X \in \C$ denote the vector
tangent to $t \mapsto p\exp(tX)$ by $\iota_p(X)$. Then the
tangent space to the fibres of the projection at $(p,q)$
is spanned by all vectors of the form $(\iota_p(X), -\iota_q(X))$.
Apply $\hat\delta(\mu)$ to the vector $(\iota_p(X), -\iota_q(X))$.
Because $\mu $ is a connection form we have
$$
d_0^*(\mu)(\iota_p(X), -\iota_q(X)) = X
$$
and
$$
d_2^*(\mu)(\iota_p(X), -\iota_q(X)) = -X.
$$
Consider the multiplication map $M \colon \hat\G \times \hat\G \to
\hat\G$. Because $\cstar$ is central we must have
$$
T_{(p_1, p_2)} M (\iota_{p_1}(X_1), \iota_{p_2}(X_2) =
\iota_{p_1p_2}(X_1 + X_2)
$$
and hence
$$
d_1^*(\mu)(\iota_p(X), -\iota_q(X)) = 0.
$$
It follows that $\hat\delta(\mu)$ is zero on vertical vectors.
As it is also clearly invariant under $\cstar$ it descends to
$\hat\G \otimes \hat\G$.  The descended form is the
tensor product connection discussed earlier and denoted
$\delta(\mu)$.

Let $\alpha = s^*(\delta(\mu))$.
We then have that
\begin{align*}
     \delta(\alpha) &= (\delta(s)^*) (\delta\delta(\mu)) \\
&= (1)^*(\delta^2(\mu))\\
&= 0
\end{align*}
as $\delta^2(\mu) $
is the flat connection on $\delta^2(P)$.   Also $d\alpha =
s^*(d\delta(\mu)) = \delta(R)$.

In more detail $\alpha$ and $R$ satisfy:
\begin{align}
d_{0}^{*}R - d_{1}^{*}R + d_{2}^{*}R &= d\alpha \label{eq:one}
\\  d_{0}^{*}\alpha - d_{1}^{*}\alpha + d_{2}^{*}\alpha -
d_{3}^{*}\alpha &= 0. \label{eq:AF}
\end{align}
Let $\Gamma(\G)$ denote the set of all pairs $(\alpha, R)$ where
$R$ is a closed, $2\pi i $ integral, two form on $\G$ and $\alpha$ is
a
one-form on $\G \times \G$ with $\delta(R) = d\alpha$ and
$\delta(\alpha) = 0$.

We have constructed
a map $C(\G) \to \Gamma(\G)$.
In the next section we construct an inverse to this
map by showing how to define a central extension from a pair
$(\alpha, R)$.  For now notice that isomorphic central
extensions with connection clearly give rise to the
same $(\alpha, R)$ and that if we vary the connection, which
is only possible by adding on the pull-back of a
one-form $\eta$ from $\G$, then we change $(\alpha, R)$
to $(\alpha + \delta(\eta), R + d\eta)$.

\subsection{Constructing the central extension}
Recall that given $R$ we can find a principal $\cstar$ bundle $P \to
\G$  with connection $\mu$ and curvature $R$ which is
unique up to isomorphism.
It is a standard result in the theory of bundles that
if $P \to X$ is a bundle with connection $\mu$ which
is flat and $\pi_1(X) = 0$ then $P$ has a section $s \colon
X \to P$ such that $s^*(\mu) = 0$. Such a section is not
unique of course it can be multiplied by a (constant) element of
$\cstar$.  Consider now our pair $(R, \alpha)$ and the bundle $P$.
As $\delta(R) = d\alpha$ we have that the connection $\delta(w) -
\pi^*(\alpha)$
on $\delta(P) \to G\times G$ is flat and hence we can find a section
$s$  such that
$s^*(\delta(w)) = \alpha$.

The section $s$ defines a multiplication by
$$
s(p, q) = p \otimes M(p, q)^* \otimes q.
$$
Consider now $\delta(s)$ this satisfies
$\delta(s)^*(\delta(\delta(w))) = \delta(s^*(\delta(w)) =
\delta(\alpha) =
0$.
  On the other hand the canonical section $1$ of $\delta(\delta(P))$
also satisfies this so they differ by a constant element of the
group.
This means that there is a $w \in \cstar$ such that for any $p$, $q$
and $r$ we must have
$$
M(M(p, q), r) = w M(p, M(q, r)).
$$
Choose $p \in \hat\G_e$ where $e$ is the identity in $\G$. Then
$M(p, p ) \in \hat\G_e$ and hence $M(p, p) = pz$ for some $z \in
\cstar$.  Now let $p=q=r$ and it is clear that we must have $w=1$.

So from $(\alpha, R)$ we have constructed $P$ and a section $s$  of
$\delta(P)$ with $\delta(s) = 1$.  However $s$ is not unique
but this is not a problem.
If we change $s$  to $s' = sz$ for
some constant $z \in \cstar$ then we have changed $M$
to $M' = M z$. As $\cstar$ is central multiplying
by $z$ is an isomorphism of central extensions with
connection. So the ambiguity in $s$  does not change the
isomorphism class of the central extension with connection. Hence
we have constructed a map
$$
\Gamma(\G) \to C(\G)
$$
as required. That it is the inverse of the earlier
map follows from the definition of $\alpha$ as $s^*(\delta(\mu))$
and the fact that the connection on $P$ is chosen so its
curvature is $R$.

\subsection{An explicit construction}

First we show how the pair $(\alpha, R)$ can be used to recover
the original central extension with connection. We will then show this
gives a construction of a central extension from any pair
$(\alpha, R) \in \Gamma(\G)$. Finally we have to show that this
gives an inverse to the map defined in the preceding section.

Let $P\G$ denote the space of all paths in $\G$ which begin at the
identity and define $\pi \colon P\G \to \G$ to be the map which
evaluates the
endpoint of the path.  We can use this map to pullback the
central extension $\hat\G \to \G$ to a central extension of
$P\G$ by $\cstar$ defined by
$$
p^{-1}(\hat\G)  = \{(f, \hat g ) \mid f(1) = \pi(\hat g) \}.
$$
Because $P\G$ is contractible this must be trivial and indeed we can
map $(f, z ) \in P\G \times \cstar$ to $(f, \hat f(1)z)$
where $\hat f$ is the (unique) lift of $f$ to a horizontal
path in $\hat\G$ starting at the identity.  The two projection
maps define a commuting diagram:
$$
\begin{array}{ccc}
     P\G \times \cstar   & \to   & \hat\G  \\
     \downarrow   &      & \downarrow \\
     P\G         & \to    & \G
    \end{array}
  $$

The product on $p^{-1}(\hat\G)$ induces a
product on $P\G \times \cstar$ which must take the form
$(f, z)(g, w) = (fg, c(f, g) zw)$ for some $\cstar$ valued
cocycle on $P\G$ which we now calculate.
Let $f$ and $g$ be
paths in $\G$ starting at the identity and  let $\hat f$ and
$\hat g$ be horizontal lifts to $\hat G$  beginning at the
identity. Then $c(f, g)$ satisfies
$$
\widehat{(fg)}(1) = c(f,g) \hat f(1) \hat g(1).
$$

Notice that  $\hat f  \hat g$ and $\widehat{fg}$
are both lifts of $fg$. Hence there is a map $\xi \colon [0,
1] \to \C^\times$  with $\xi(0) = 1$ and $\xi(1) = c(f,g)$ such that
$\hat f \hat g = \widehat{fg}\xi$.  As $\widehat{fg}$  is horizontal
it is straightforward to integrate and prove that
$$
c(f, g) = \exp(-\int_{\hat f\hat g} \mu).
$$
Let $(f, g)$ denote the path in $\G \times \G$. From the definition
of $\alpha$ we see that
\begin{align*}
     \int_{(f,g)} \alpha &=  \int_{(\hat f,\hat g)} \pi^*(\alpha) \\
                   &= \int_{\hat f} \mu - \int_{\hat f\hat g} \mu
		                          + \int_{\hat g} \mu
		\end{align*}
Using the fact that $\hat f$ and $\hat g$ are horizontal we find that
$$
c(f, g) = \exp(\int_{(f,g)} \alpha).
$$

As in \cite{Mur} we can now identify the kernel of the homomorphism
$ P\G \times \C^\times \to \hat G$.  This is all pairs $(h, z)$ such
that the holonomy of the connection is equal to $z^{-1}$.  As we are
assuming that $\G$ is simply connected we can extend any loop $h$ to
a
map $\tilde h $ from  a  disk $D$ into $\G$ and define
$$
H(h, R) = \exp(\int_{\tilde h(D)} R ).
$$
Notice that the $2\pi i$ integrality of $R$ implies that
$H(h, R)$ is well-defined.  The kernel is therefore the
subgroup of all pairs $(h, H(h, R)^{-1})$.

We can now see how to define a central extension given the pair
$(\alpha,
R)$.  First we define $c(f, g)$ by
$c(f, g) = \exp(\int_{(f,g)} \alpha)$ and it can be checked, as in
\cite{Mur}, that this co-cycle makes $P\G \times \C^\times$ into a
group.  Then define $H(h, R)$ as above and consider the subset of
all pairs $(h, H(h, R)^{-1})$. It can be shown that this is a normal
subgroup and the quotient defines the central extension.

We define a connection 1-form $\mu$ on the principal
$\cstar$ bundle $\hat{\G}$ using the same
technique as in \cite{Mur1}.  We use the map
$$
P\G \times \cstar \to \hat\G
$$
to pullback the connection one-form $\mu$ on $\hat G$ to a
connection one-form $\hat\mu$.
A straightforward calculation shows that this is given by
$$
  z^{-1}dz + \hat \mu
$$
where
$$
\hat{\mu}  = \int_{[0,2\pi]}\iota_
{(0,\frac{d}{dt})}\ev^{*}R
$$
and  $\ev\colon \P\G\times [0,2\pi]\to \G$ is
the evaluation map $\ev(f,t) = f(t)$.

In the case that we start with a pair $(\alpha, R)$ it is
straightforward
to check that this connection descends to a connection on $\hat \G$.

We now have a procedure for constructing from a given central
extension a pair $(\alpha, R)$ and from a pair $(\alpha, R)$ a
central extension.  It is clear from the construction that
if we start with a central extension, construct
$(\alpha, R)$ and then construct a central extension we get back to
were we started from.  Consider now what happens if we
start with an $(\alpha, R)$, construct the central extension and
then construct a pair $(\tilde \alpha, \tilde R)$.  It is
straightforward
to  show that we have
\begin{align*}
     \pi^*(\tilde \alpha) &= c^{-1} dc + \delta(\hat\mu) \\
    \hbox{and} \quad\quad  \pi^*(\tilde R)&= d\hat\mu.
\end{align*}
Using the definition of $c$ and $\hat \mu$
we can show that $\alpha$ and $R$ satisfy the same equations and
hence
deduce that $\tilde \alpha = \alpha $ and $\tilde R = R$ as $\pi^*$
is injective.

\section{Loop groups}
There are a number of variants of the loop group that we wish to
consider.  To define these let $K$ be a compact group and
consider first
$$
L(K) = \{ \gamma \colon [0, 2\pi] \to K \mid \gamma(0)  = \gamma(2\pi) \}
$$
this has a subgroup of based loops
$$
L_0(K) = \{ \gamma \colon [0, 2\pi] \to K \mid \gamma(0)  = \gamma(2\pi)
= 1 \}.
$$
We assume here that the maps $\gamma$ are smooth on $[0, 2\pi]$.
Now map $[0, 2\pi]$ to the circle $S^1$ by $\theta \to
\exp(i\theta)$.  We therefore have
$$
\Omega(K) = C^\infty(S^1, K)
$$
the space of smooth maps from the circle to $K$ with a subgroup
of $LK$ and we let $\Omega_0(K) = \Omega(K) \cap L_0K$.   These
groups
are all Frechet Lie groups and their Lie algebras are the analogous
spaces of maps of $[0, 2\pi]$ to $\mathfrak{k}$ and denoted
by  $L(\mathfrak{k})$, $L_0(\mathfrak{k})$,
$\Omega(\mathfrak{k})$ and $\Omega_0(\mathfrak{k})$.

In the case where $\G = L (K)$ there is a
well known expression for the curvature $R$ of a
left invariant connection on $\hat{L(K)}$
--- see \cite{PreSeg}.  We can also write down
a 1-form $\alpha$ on $L(K)\times L(K)$
such that $\d(R) = d\alpha$ and $\d(\alpha) = 0$. These are:
\begin{align}
     R &= \frac{i}{4\pi}\int_{S^1} \< \Theta, \partial_\theta
\Theta\>
     d\theta \\
     \alpha &= \frac{i}{2\pi} \int_{S^1} \< d_2^*\Theta, d_0^*Z \>
d\theta
\label{eq:newRalpha}
\end{align}

Here $\Theta$ denotes the
Maurer-Cartan form on the Lie group $K$, that is
$\Theta(k)( kX) = X$,  $Z $ is the function on $LG$ defined by $Z(g) =
(\partial_\theta
g ) g^{-1}$ and  $\<\ ,\ \>$ is an
  invariant inner product normalised
so that the  longest root has length squared equal to $2$.  Note that 
$F$ is left invariant and that $A$
is left invariant in the first factor of  $\G\times \G$.

The formulae \eqref{eq:newRalpha} contain implicit wedge products and to
avoid confusion let us explain what these are.  Let $\omega_i$ be
differential forms of degree $d_i$ with values in a vector space
$V_i$ for $i=1, \dots, k$.
If $p \colon V_1 \times \dots \times V_k \to \C$
is a $k$ linear map then we define
$p(\omega_1, \dots, \omega_k)$,
a differential form of degree $d = d_1 + \cdots + d_k$
by:\begin{multline}
p(\omega_1, \dots, \omega_k)(X_1, \dots, X_k) \\
= \sum_{\pi \in S_d} \text{sign}({\pi}) p(\omega_1(X_{\pi(1)}, \dots, 
X_{\pi(d_1)}) ,
\dots, \omega_k(X_{\pi(d-d_k+2)}, \dots, X_{\pi(d)}))
\end{multline}
where the sum is
over all permutations $\pi$ of $\{1, \dots, d\}$
with $\text{sign}(\pi)$ the sign of the permutation.
Note the potential confusion that if each of the $\omega_i$ is
of degree $1$ and $p$ is an
antisymmetric function then $p(\omega_1(X_1), \dots, \omega_d(X_d))$ 
is already an
anti-symmetric function of $X_1, \dots, X_d$. It is, in fact,
$1/d!$ times $p(\omega_1, \dots, \omega_d)$
applied to $X_1, \dots, X_d$. For later  use we record here some 
identities relating
  $\Theta$ and $Z$.
At a point $g \in LK$ we have
\begin{equation}
     \partial_\theta \Theta = \ad(g^{-1}(dZ))
\label{eq:identone}
\end{equation}
  and if $X$ is in $L\mathfrak{k}$ then
\begin{equation}
  \partial_\theta(\ad(g^{-1}(X)) = \ad(g^{-1})([X, Z]) + \ad(g^{-1})
  \partial_\theta X.
\label{eq:identtwo}
\end{equation}

\section{Lifting bundle gerbes and the string class}
\label{sec:lifting}

\subsection{The string class}

Consider a central extension
$$
\cstar \to \hat{\G}\stackrel{\pi}{\to} \G
$$
and let $P$ be a principal bundle over a manifold $M$
with structure
group $\G$. There is a characteristic
class in $H^3(M, \Z)$ which is
the obstruction to lifting $P$
to a $\hat\G$ bundle. This is easily
described in Cech cohomology.
  Let $\{U_\a\}_{\alpha \in I}$ be a
good
cover of $M$ with respect to which $P$ has  transition
functions $g_{\a\b}$. As the double intersections are contractible we
can lift the
transition functions to  $\hat g_{\a\b}$ with values in
$\hat\G$ such
that $\pi \hat g_{\a\b} = g_{\a\b}$. Because
$g_{\b\c}g^{-1}_{\a\c}g_{\a\b} = 1$ it follows that
$$
d_{\a\b\c} =
\hat g_{\b\c} \hat g^{-1}_{\a\c} \hat g_{\a\b}
$$
which defines a class in $H^2(M, \Cstar)$ and it is straightforward
to check that
this class is zero if and only if the bundle lifts to
$\hat \G$ \cite{Mur}.  A
standard calculation with the short exact sequence
$$
\Z  \to \C \to \Cstar
$$
identifies this class with a class in
$H^3(M, \Z)$.  In that
case that $\G$ is the loop group and $\hat
\G$ the central
extension of the loop group this class is called the
{\em string class}
\cite{Kil}.

\subsection{Lifting bundle gerbes}

We first briefly review the
notion of a
\emph{bundle gerbe} from \cite{Mur1} cast in the
language of simplicial line bundles. Let
$\pi:Y\to M$ be a local split map. Recall that
it defines a simplicial space $Y^{[p]}$.  A bundle
gerbe over  a manifold $M$ is a pair $(P, Y)$
where $P$ is a simplicial line bundle over the simplicial
space defined by $Y$. This means that $P$ is a
$\cstar$ bundle over $Y^{[2]}$.

A bundle gerbe $P$ has a characteristic class
$DD(P)$ in $H^{3}(M;\Z)$ associated to it --- the Dixmier-Douady
class.  We refer to \cite{Mur1} for the definition
and properties of $DD(P)$.
We can realise the image of the Dixmier-Douady
class in real cohomology in a manner analogous to the
chern class.  A bundle gerbe $(P, Y, M)$ can be equipped with a
\emph{bundle gerbe connection}.  This is a
connection $\nabla$ on $P$ which is compatible
with the bundle gerbe product in the sense that $\delta(\nabla)$
is the flat connection for the trivialisation $s$ of $\delta(P)$.
It is not hard to see that bundle gerbe
connections always exist, indeed if $\nabla$ is any
connection then $\delta(s^*(\nabla)) =  0$
and hence from \eqref{eq:complex} $s^*(\nabla) = \delta(\rho)$ for
some $\rho $
on $Y^{[2]}$ and $\nabla - \pi^*(\rho)$ is a bundle gerbe
connection.   If
$F_{\nabla}$ denotes the curvature of the
bundle gerbe connection $\nabla$ then it is
easy to see that we have $\d(F_{\nabla}) = 0$.
Hence we can solve the equation $F_{\nabla} =
\d(f)$ for some $f\in \Omega^{2}(Y)$.  A choice
of $f$ is called a \emph{curving} for the bundle gerbe
connection $\nabla$.  Since $F_{\nabla}$ is
closed, we get $\d(df) = 0$ and so we have
$df = 2\pi i\pi^{*}\omega$ for some necessarily
closed 3-form $\omega$ on $M$.  In \cite{Mur1}
it is shown that $\omega$ is an integral 3-form
and represents the image of the Dixmier-Douady
class $DD(P)$ of $\omega$ in $H^{3}(M, \R)$.

Suppose that $G$ is a Lie group forming part of a
central extension
$$
\cstar \to \hat{\G} \to \G.
$$
Recall from \cite{Mur1} that we can associate to a
principal $\G$ bundle $P(M,\G)$
a bundle gerbe $(\tilde{P},P,M)$ --- the so called
\emph{lifting bundle gerbe}.  $\tilde{P}\to P^{[2]}$
is the pullback $\tilde{P} = \tau^{-1}\hat{\G}$ of
$\hat{\G}\to \G$ by the natural map $\tau:P^{[2]}\to \G$
defined by $p_{2} = p_{1}\tau(p_{1},p_{2})$ for $p_{1}$
and $p_{2}$ points of $P$ lying in the same fibre.
$\tau$ satisfies the property $\tau(p_{1},p_{2})
\tau(p_{2},p_{3}) = \tau(p_{1},p_{3})$ for points
$p_{1}$, $p_{2}$ and $p_{3}$ all lying in the same fibre.
$\tilde{P}$ inherits a bundle gerbe product from the
product in $\hat{G}$.  The Dixmier-Douady class of the
bundle gerbe $\tilde{P}$ measures the obstruction to
lifting the structure group of the principal $G$ bundle
to $\hat{\G}$.  It follows that when  $\G$ is a loop
group the
Dixmier-Douady class is (the image in real cohomology of) the string
class of the bundle  $P$.

 From the principal bundle we can construct a simplicial
space as above and if
$$
\tau \colon P^{[2]} \to \G
$$
is defined by $p_2= p_1 \tau(p_1, p_2)$
then we can define
$$
\tau  \colon P^{[k+1]} \to \G^{k}
$$
by
$$
\tau(p_1, \dots, p_{k+1}) = (\tau(p_1, p_2), \dots, \tau(p_{k},
p_{k+1}).
$$
It is straightforward to check that this is a {\em simplicial}
map, that is it commutes with the face and degeneracy maps.  It
follows
that pullback of differential forms by $\tau^*$ commutes with
$\delta$.
Suppose that $\mu$ is a connection on $\hat{G}$.  Then
the natural connection $\tilde{\mu} =  \tau^{*}\mu$ on $\tilde{P}$
is not a bundle gerbe connection on $\tilde{P}$.
Indeed  from the equation $\delta(\mu) =\pi^*(\alpha)$
we see that $\delta(\tau^{*}(\mu)) =
\tau^*(\alpha)$. However the form $\beta = \tau^*(\alpha) $
satisfies $\delta(\beta) = \tau^*(\delta(\alpha)) = 0$
as $\delta(\alpha) = 0$. So we can solve  the
equation $\d(\epsilon) = \beta$ for some 1-form $\epsilon$ on
$P^{[2]}$.
Then the connection $\tilde{\mu} -
\epsilon$ is a bundle gerbe connection on $\tilde{P}$. Its
curvature is given by  $\tau^*(R) - d\epsilon$, where $R = d\mu$ is
the  curvature of $\mu$. The curving is therefore
a two-form $f$ on $P$ satisfying
  $\d(f) = \tau^*(R) - d\epsilon$
for some 2-form $f$ on $P$.
 From $f$ the Dixmier-Douady class $\omega$ is obtained
as $df = \pi^*(\omega)$.

To proceed further we need to concentrate on  a specific example
so we will let $\G = L(K)$.  So we  have a principal $L (K)$ bundle
$P(M,L(K))$ on $M$ for $K$ a compact Lie group, we can form the
lifting bundle gerbe $(\tilde{P},P,M)$ associated to
the central extension
$$
\cstar \to \hat{L} (K) \to L (K)
$$
where $\hat{L} (K)$ is the Kac-Moody group.  We can
form a bundle gerbe connection $\nabla$ on $\tilde{P}$
in the manner described above using the natural connection
on $\hat{L} (K)$ and the 1-form $\alpha$ on $L (K)
\times L (K)$.  We will show that it is possible
to write down an expression for the three curvature $\omega$
of the bundle gerbe connection $\nabla$.

Suppose we have chosen a connection 1-form $A$ on the
principal $\G$ bundle $P\to M$.  Then this is a one-form on $P$ with
values in $\mathfrak{g}$.  It is straightforward to show that
\begin{equation}
     \pi_1^*(A) = \ad(\tau^{-1})\pi_2^*(A) + \tau^*(\Theta).
\label{eq:Atau}
\end{equation}
Let $\tau_{ij}(p_1, \dots, p_k) = \tau(p_i, p_j)$. Then using
$\beta = (\tau_{12}\times \tau_{23})^{*}\alpha$, the definition
of $\alpha$ from \eqref{eq:newRalpha} and the identity \eqref{eq:Atau}
we obtain
\begin{equation*}
\beta
= \frac{i}{2\pi} \int_{S^{1}} \< \pi_{13}^*A -
\ad(\tau_{12})^{-1}\pi_{23}^*A ,
\partial_\theta (\tau_{23} )\tau_{23}^{-1} \> d\theta
\end{equation*}
where $\pi_{12}(p_{1},p_{2},p_{3}) = p_3$, etc.

Define a 1-form  $\epsilon$ on $P^{[2]}$ by
\begin{equation}
     \epsilon = \frac{i}{2\pi}\int_{S^1} \< \pi_2^*A, \tau^*(Z)\>
d\theta.
\label{eq:epsilon}
\end{equation}
Using the fact that $\tau_{12}\tau_{23} = \tau_{13}$ we obtain
$$
\tau_{12}^*(Z) +
\ad(\tau_{12})(\tau_{23}^*(Z)) =
\tau_{13}^*(Z)
$$
and with this and \eqref{eq:Atau} we can show that $\delta(\epsilon)
= \beta$.

To solve the equation $\tau^{*}R - d\epsilon = \d(f)$ for some
choice of curving $f$ we first need an explicit expression for
$\tau^{*}R - d\epsilon$.  Using \eqref{eq:newRalpha} we have that
\begin{equation*}
     \tau^*R = \frac{i}{4\pi}\int_{S^1} \<\tau^*\Theta,
     \partial_\theta(\tau^*\Theta) \>  d\theta
     \end{equation*}

Using the standard fact that $d\Theta = - (1/2) [ d\Theta, d\Theta]$ and the
identities \eqref{eq:identone}
and
\eqref{eq:identtwo} we can show that we have
\begin{multline*}
\tau^*{R} - d\epsilon
  =     \frac{i}{4\pi} \int_{S^1}  \< \pi_1^*(A), \partial_\theta
\pi_1^*(A)\> - \< \pi_2^*(A),
\partial_\theta \pi_2^*(A)\>  \\
  -\<[\pi_2^*(A), \pi_2^*(A)], \tau^*(Z)\> - 2\< \pi_2^*(dA),
\tau^*(Z) \>
d\theta.
\end{multline*}

Recalling that $F$, the curvature of $A$ satisfies $F = dA + 1/2[A,
A]$
we have that $\tau^*{R} - d\epsilon $ is equal to
\begin{equation*}
\delta\left(\frac{i}{4\pi} \int_{S^1}
\< A, \partial_\theta A \>\right)
- \frac{i}{2\pi}\int_{S^1} \<\pi^*_2(F) , \tau^*(Z) \> d\theta
\end{equation*}
where $\delta = \pi_1^* - \pi^*_2$. We want to solve $\tau^*{R} -
d\epsilon  =
\delta(f)$ so we now need to  write the two-form
$$
\frac{i}{2\pi}\int_{S^1} \<\pi^*_2(F) , \tau^*(Z) \> d\theta
$$
as $\d$ of a
two-form on $P$.  To this end, choose  a map $\Phi \colon
P\to C^\infty([0, 2\pi], \mathfrak{k})$ satisfying
\begin{equation}
\Phi(pg) = \ad({g^{-1}})\Phi(p) + g^{-1}\frac{\partial g}
{\partial \theta}.   \label{eq:thirtyfour}
\end{equation}
We call any such function a {\em twisted Higgs field}
for $P$. As a convex combination of twisted  Higgs fields
is a twisted Higgs field and they clearly exist when $P$
is trivial it is straightforward to use a partition
of unity on $M$ to construct a twisted Higgs field for $P$.
Notice that a twisted Higgs field will not generally
live in $L(\mathfrak{k})$. However if we start with a bundle
$P(M, \Omega(K))$ then the twisted Higgs field can be required to
take values in $\Omega(\mathfrak{k})$.

Choose then $\Phi$ a twisted Higgs field for $P$. It
satisfies
$$
\ad(\tau)\pi^*_1(\Phi) = \pi_2^*(\Phi) +  \tau^*(Z)
$$
and hence
\begin{align*}
     \<\pi^*_2(F) , \tau^*(Z) \> &= \<\pi^*_2(F) ,
     \ad(\tau)\pi^*_1(\Phi) - \pi_2^*(\Phi) \> \\
     &=  \<\ad(\tau^{-1})\pi^*_2(F) ,  \pi^*_1(\Phi)\>
     - \< \pi^*_2(F), \pi_2^*(\Phi) \>\\
     & = \pi_1^*(\<F, \Phi\>) - \pi_2^*(\<F, \Phi\>).
  \end{align*}

  Finally
we can write $\tau^*R - d\epsilon = \d(f)$ where $f$
is the two-form on $P$ defined by
\begin{equation*}
f = \frac{i}{2\pi} \int_{S^1} (\frac{1}{2}\<A, \partial_\theta A \>
- \<F, \Phi\> ) d\theta
\end{equation*}
and hence
\begin{equation*}
     df = \frac{i}{2\pi} \int_{S^1}  \<dA, \partial_\theta A\> - \<
     dF, \Phi \> - \<F, d\Phi\> d\theta
     \end{equation*}
Using the Bianchi identity for $F$ we obtain
$$
df = -\frac{i}{2\pi} \int_{S^1} \<F, \nabla \Phi\> d\theta
$$
where
$$
\nabla \Phi = d\Phi + [A, \Phi] - \partial_\theta A.
$$
It is straightforward to show that
$$
\pi^*_1(\nabla\Phi) = \ad(\tau^{-1}) \pi^*_2(\nabla\Phi)
$$
which shows that $\nabla\Phi$ descends to a one form on $M$
with values in the adjoint bundle of $P$, just as $F$ descends
to a two-form on $M$ with values in the adjoint bundle
of $P$.

Finally we have
\begin{theorem}
     Let $P \to M$ be a principal $L K$ bundle and let $A$
     be a connection for $P$ with curvature $F$ and  $\Phi$
     be a twisted Higgs field for $P$. Then the
     string class of $P$ is represented in de Rham cohomology
     by the three-form
     \begin{equation}
\label{eq:stringclass}
  -\frac{1}{4\pi^2} \int_{S^1} \<F, \nabla \Phi\> d\theta
\end{equation}
where
\begin{equation}
\label{eq:nablaHiggs}
\nabla \Phi = d\Phi + [A, \Phi] - \partial_\theta A.
\end{equation}\end{theorem}

\subsection{Reduced splittings}
\label{sec:gomi}
In \cite{Gom} the concept of {\em reduced splitting} is used to
calculate a connection and curving for the lifting bundle
gerbe.  In this section we explain how this fits in with
our constructions in the case of the loop group.

First we have
\begin{definition}[\cite{Gom}]
Define the {\em group cocycle} $Z \colon L(K) \times L(\LK) \to i\R$ by
$$
(0, Z(g, X)) = \ad(\hat g^{-1})((X,0))- (\ad(g^{-1}) (X ) , 0)
$$
\end{definition}
where here we assume a splitting of $\widehat{L(K)} $
as $L(K) \oplus i\R$ and $\hat g$ is a lift of $g$ to the
central extension. The group cocycle  $Z$ encapsulates the information
of the central extension in a similar manner to our $\alpha$. Indeed in
the case of the loop group  with the extension defined by the $(R, \alpha)$
in \eqref{eq:newRalpha} then we have
$$
Z(g^{-1}, X) = - \alpha(1,g)(X,0) = \frac{i}{2\pi} \int_{S^1} \< X, 
\partial_\theta (g) g^{-1}\> d\theta.
$$

\begin{definition}[\cite{Gom}]
A {\em reduced splitting} for a loop group principal bundle 
$P(M,L(K))$ is a map
$$
\ell \colon P \times L(\LK) \to \R
$$
which is linear in the second factor and satisfies $\ell(p, X) = 
\ell(pg, \ad(g^{-1})(X)) +
Z(g^{-1}, X)$.
\end{definition}

A straightforward calculation shows that if $\Phi$ is a Higgs field then
$$
\ell(p, X) =  \frac{i}{2\pi} \int_{S^1} \< \Phi, X\>d\theta
$$
is a reduced splitting.

\subsection{The path fibration}
We will illustrate the above discussion for the
case of the $L_0(K)$ bundle $\P K(K,L_0(K))$,
where $\P K\to K$ is the path fibration of the group
$K$ --- see \cite{CarMur}.  In this case there is a
canonical closed 3-form on $K$ representing the
image of the generator of $H^{3}(K;\Z) = \Z$ inside
$H^{3}(K;\Reals)$ given by
$$
\omega_3 = \frac{1}{48\pi^2} \<[\widehat\Theta, \widehat\Theta], 
\widehat\Theta\>
$$
where $\widehat\Theta$ is the {\em right} invariant
Maurer-Cartan form on $K$
see for example \cite{Bry}.  We will show that for the
path fibration the string class \eqref{eq:stringclass} is $\omega_3$.

First we need to define a connection.  The tangent
space at $p \in \P K$ is the set of all vector
fields along the path $p$ and  such a vector field
is vertical if it vanishes at $2\pi$. A right invariant
splitting is given by
$$
H_{p} = \{\theta \mapsto (\theta/2\pi)R_{p(\theta)}X \mid  \in \LK\}.
$$
Clearly this satisfies $R_{g_{*}}
H_{p} = H_{pg}$ and the corresponding connection
one-form is
$$
A = \Theta - \frac{\theta}{2\pi} \ad(p^{-1}) \pi^*(\widehat\Theta).
$$
A calculation shows that the curvature of $A$ is
$$
F = \left( \frac{\theta^2}{8 \pi^2} - \frac{\theta}{4 \pi}\right)
\ad(p^{-1}) [\pi^*(\widehat\Theta), \pi^*(\widehat\Theta)].
$$
A suitable Higgs field is given by
$$
\Phi(p) = p^{-1}\frac{\partial p} {\partial \theta}.
$$
and another  calculation shows that
$$
\nabla\Phi = \frac{1}{2\pi} \ad(p^{-1}) \pi^*(\widehat\Theta).
$$
Putting these into the formula for the string class 
\eqref{eq:stringclass} gives the  required result.

\section{Calorons and the string class}
In \cite{GarMur} as part of the study of calorons a correspondence was
introduced between $K$ bundles $\tilde P$ on a manifold $M \times
S^1$ and
$\Omega(K)$ bundles $P$ on $M$.  This correspondence also related a
connection $\tilde A$ on $\tilde P$ to a connection $A$ on $P$ and
a section $\Phi$ of the twisted adjoint bundle. We will describe this
correspondence and show that the integral over the circle of the
Pontrjagin class of $\tilde A$ is the representative for the
string class \eqref{eq:stringclass}.

If $P \to M$ is an $\Omega(K)$ bundle we define
\begin{equation*}
     \tilde P = (P \times K \times S^1)/\Omega(K)
    \end{equation*}
    where $\Omega(K)$ acts on $P\times K \times S^1$ by
    $g(p, k, \theta) = (pg^{-1}, g(\theta)k, \theta)$.
   Letting $[p, k, \theta]$ denote the equivalence class
   of $(p, k, \theta)$ we have a right $K$ action on $\tilde P$
   given by $[p, k, \theta]h = [p, kh, \theta]$ and
   a projection map $\pi \colon \tilde P \to M \times S^1$ defined by
   $\pi([p, k, \theta])= (\pi(p), \theta)$. The fibres of $\pi$ are
the
   orbits of $K$  and $\tilde P$ is a principal $K$ bundle.

Starting instead with a principal $K$ bundle $\tilde P \to M \times
S^1$ we define a bundle $P \to M$ whose fibre at $m$ is all the
sections of $\tilde P$ restricted to $\{m\} \times S^1$.  This is
clearly acted on by $\Omega(K)$.

Given a connection 1-form $A$ on $P(M,\Omega(K))$ we can
define a connection 1-form $\tilde{A}$ on $\tilde{P}$ by
pushing forward the 1-form on $P\times K\times S^{1}$
which is also denoted by $\tilde{A}$ and is
given by
\begin{equation}
\tilde{A} = \ad(k^{-1})(A)+ \Theta + \ad(k^{-1}){d\theta\Phi} . 
\label{eq:fourtythree}
\end{equation}
To check that we can push forward $\tilde{A}$
to a connection 1-form on $\tilde{P}$ we first have to
show that $\tilde{A}$ is invariant under the
action of $L (K)$ on $P\times K\times S^{1}$ and
that $\tilde{A}$ kills vectors tangent to the
fibering $P\times K\times S^{1}\to \tilde{P}$.  Let
$(p,k,\phi)$ be a point of $P\times K\times S^{1}$
and suppose that $(X,\eta,\lambda)$ is a tangent
vector at $(p,k,\phi)$.  Then if $g\in \Omega(K)$
we have
\begin{equation}
R_{g}^{*} \tilde{A}((p,k,\phi);(X,\eta,\lambda))
=  \ad(k^{-1})(A(p;X)(\phi)) + L_{k^{-1}}(\eta) +
\ad(k^{-1})(\lambda \Phi(\phi)).
\end{equation}
so $\tilde{A}$ is invariant under the action
of $L (K)$.  We now check that $\tilde{A}$
kills vectors tangent to the fibering $P\times K\times
S^{1}\to \tilde{P}$.  The vectors tangent to this
fibering are of the form $(\iota_{p}(X), - X(\theta)k,0)$
where $X $ is in $\Omega(\LK)$ the Lie algebra of $\Omega(K)$ and
where $\iota_{p}(X) = \frac{d}{dt}|_{t=0}p\exp(tX)$.  It
is easy to see that these vectors lie in the kernel
of $\tilde{A}$.  We now need to check that the
pushed forward 1-form on $\tilde{P}$ (which we will
also denote by $\tilde{A}$) is a connection
1-form on $\tilde{P}$.  It is a straightforward matter
to check that $\tilde{A}$ is equivariant under
the action of $K$ on $\tilde{P}$.  Now let $\xi \in \LK$
and let $\tilde{p} = (p,k,\phi)$ be a point of $P\times
K\times S^{1}$.  Then we have $\iota_{\tilde{p}}(\xi) =
(0,L_k(\xi),0)$ and $\tilde{A}(0,L_k(\xi),0) = \xi$.
Hence $\tilde{A}$ pushes forward to define a
connection 1-form on $\tilde{P}$.

If we compute the  curvature $\tilde{R} =
d\tilde{A} + 1/2[\tilde{A},\tilde{A}]$
for the connection $\tilde{A}$ we obtain:
\begin{equation}
\tilde R = \ad(k^{-1})(F + \nabla(\Phi)d\theta)
\end{equation}
where $\nabla(\Phi)$ is defined as in \eqref{eq:nablaHiggs} and
$F$ denotes the curvature of the connection
$A$ on the principal $L (K)$ bundle $P$.
The Pontrjagin form of the connection $A$ on the
bundle over $S^1 \times M$ is
\begin{equation}
-\frac{1}{8 \pi^2} \<\tilde R, \tilde R\> =
-\frac{1}{8\pi^2}\left(\<F, F\> + 2\< F, \nabla\Phi\>\right)
\label{eq:fourtyfive}
\end{equation}
and integrating over the circle gives the
string class of the bundle $P(M, \Omega(K))$.
We have hence proved

\begin{theorem} Let $P \to M$ be an $\Omega(K)$ bundle
     and $\tilde P \to M \times S^1$ the corresponding $K$
     bundle. Then the string class of $P$ is obtained from
     the Pontrjagin class of $\tilde P$ by integrating
     over the circle.
     \end{theorem}

This result enables an easy proof of a result of Killingback
also proved in \cite{CarMur, Gom}.  We start with a $K$
bundle $Q \to X$ and let $P$ and $M$ be the loop spaces
of $Q$ and $X$ respectively. Then $P$ is an $\Omega(K)$
bundle over $M$ and Killingback shows that the
string class of $P$ is obtained from the Pontrjagin
class of $Q$ by pulling back by the evaluation
map
\begin{equation}
     \ev \colon M \times S^1 \to X
\label{eq:evmap}
\end{equation}
and integrating over the circle.   If $P$ is the space
of loops in $Q$ there is a natural map
$ P \times K \times S^1 \to Q $
given by $(p, k, \theta) \mapsto p(\theta)k$.  This is constant
on the orbits of the $\Omega(K)$ action and hence defines a
map $\tilde P \to Q$ which is $K$ equivariant and covers the
evaluation map \eqref{eq:evmap}. This shows that $\tilde P $ is the
pull-back of $Q$ by the evaluation map. Hence the Pontrjagin class
of $\tilde P$ is the pull-back by the evaluation map of the
Pontrjagin class of $Q$. Thus the string class of $P$ is the
integral over the circle of the pull-back by the evaluation
map of  the Pontrjagin class of $Q$.

\section{Conclusion}
If we try and apply the caloron construction to an $L(K)$ bundle
with a connection we obtain a $K$ bundle on $S^1$ which is only
smooth on $(0, 2\pi)$.  There should be a theory of $K$ bundles
on $[0, 2\pi]$ with connection which patch together over $\{0\} \times M$
and $\{2\pi\} \times M$ in such as way as to recover the result of
\cite{GarMur} relating $L(K)$  bundles on $M$ to such bundles.

Notice that this approach could be used to show that any $L(K)$ bundle
over $M$ has a
three class  which we could define by transgressing the Pontrjagin class
of the induced $K$ bundle over $S^1 \times M$.  However it would not be
clear that this three class was the string class, the obstruction to
lifting to $\widehat L(K)$. For this we needed the theory of  bundle 
gerbes  and the
connections and curvings which provide a bridge between the Cech 
description of the
string class and a de Rham realisation of it.

\end{document}